\documentclass[11pt]{amsart}
\usepackage{amsfonts}
\usepackage{amsmath}
\usepackage{amssymb}
\usepackage{color}
\usepackage{enumerate}

\newtheorem{theorem}{Theorem}

\newtheorem{corollary}[theorem]{Corollary}

\newtheorem{definition}[theorem]{Definition}
\newtheorem{example}[theorem]{Example}

\newtheorem{lemma}[theorem]{Lemma}

\newtheorem{proposition}[theorem]{Proposition}

\def\a{\alpha}

\def\e{\epsilon}

\def\n{\nu}

\def\F{\mathbf F}

\def\ot{\otimes}
\def\ra{\rightarrow}

\def\QM{{\mathbf Q}}
\def\CM{{\mathbf C}}

\def\NM{{\mathbf N}}

\begin{document}

%%%%%%%%%%%%%%%%%%%%%%%%%%%%%%%%%%%%%%%%%%%%%%%%%%%%%%%%%%%%

\title[Logarithmic derivatives]{Logarithmic derivatives and generalized Dynkin operators.} 

\date{June 6, 2012}

%%%%%%%%%%%%%%%%%%%%%%%%%%%%%%%%%%%%%%%%%%%%%%%%%%%%%%%%%%%%

\author{Fr\'ed\'eric Menous}
% \address{.}
% \email{}
% \urladdr{}

\author{Fr\'ed\'eric Patras}
% \address{Laboratoire J.-A.~Dieudonn\'e
%          		UMR 6621, CNRS,
%          		Parc Valrose,
%          		06108 Nice Cedex 02, France.}
% \email{patras@math.unice.fr}
% \urladdr{www-math.unice.fr/~patras}

%%%%%%%%%%%%%%%%%%%%%%%%%%%%%%%%%%%%%%%%%%%%%%%%%%%%%%%%%%%%
\begin{abstract}
Motivated by a recent surge of interest for Dynkin operators in mathematical physics and by problems in the combinatorial theory of dynamical systems,
we propose here a systematic study of logarithmic derivatives in various contexts. In particular, we introduce and investigate generalizations of the Dynkin operator for which we obtain Magnus-type formulas.
\end{abstract}

%%%%%%%%%%%%%%%%%%%%%%%%%%%%%%%%%%%%%%%%%%%%%%%%%%%%%%%%%%%%
\maketitle
\tableofcontents
%%%%%%%%%%%%%%%%%%%%%%%%%%%%%%%%%%%%%%%%%%%%%%%%%%%%%%%%%%%%

\section*{Introduction}
\label{sect:intro}
Dynkin operators are usually defined as iterated brackettings. They are particularly popular in the framework of linear differential equations and the so-called continuous Baker-Campbell-Hausdorff problem (to compute the logarithm of an evolution operator). We refer to \cite{reutenauer} for details and an historical survey of the field.
Dynkin operators can also be expressed as
a particular type of logarithmic derivatives (see Corollary~\ref{classical} below).
They have
received increasingly more and more interest during the recent years, for various reasons.
\begin{enumerate}
 \item In the theory of free Lie algebras and noncommutative symmetric functions, it was shown that Dynkin operators generate the descent algebra (the direct sum of Solomon's algebras of type $A$, see \cite{gelfand,reutenauer}) and play a crucial role in the theory of Lie idempotents. 
 \item Generalized Dynkin operators can be defined in the context of classical Hopf algebras \cite{patreu2002}. The properties of these operators generalize the classical ones. Among others (see also \cite{EGP,EMP}), they can be used to derive fine properties of the renormalization process in perturbative quantum field theory (pQFT): the Dynkin operators can be shown to give rise to the infinitesimal generator of the differential equation satisfied by renormalized Feynman rules  \cite{EFGP}. This phenomenon has attracted the attention on logarithmic derivatives in pQFT where generalized Dynkin operators are expected to lead to a renewed understanding of Dyson-Schwinger-type equations, see e.g. \cite{KY06}.
\end{enumerate}

The present article was however originated by different problems, namely problems in the combinatorial theory of dynamical systems (often referred to as Ecalle's \it mould calculus\rm , see e.g. \cite{sauzin}, also for further references on the subject) and in particular in the theory of normal forms of vector fields. 
It appeared very soon to us that the same machinery that had been successfully used in the above mentioned fields and problems was relevant for the study of dynamical systems.
However, the particular form of logarithmic derivatives showing up in this field 
requires the generalization of the known results on logarithmic derivatives and Dynkin operators to a broader framework: we refer to the section~\ref{exemple} of the present article for more details on the subject. 

The purpose of the present article is therefore to develop further the algebraic and combinatorial theory of logarithmic derivatives, following various directions and point of views (free Lie algebras, Hopf algebras, Rota-Baxter algebras...), all of them known to be relevant for the study of dynamical systems but also to various other fields, running from the numerical analysis of differential equations to pQFT.

We limit the scope of the present article to the general theory and plan to use the results in forthcoming articles.

\section{Twisted Dynkin operators on free Lie algebras}

Let $T(X)$ be the tensor algebra over an alphabet $X=\{x_1,...,x_n,...\}$. 
It is graded by the length of words: the degree $n$ component $T_n(X)$ of $T(X)$ is the linear span (say over the rationals) of the words $y_1...y_n,\ y_i\in X$. The length, $n$, of $y_1...y_n$ is written $l(y_1...y_n)$.
It is equipped with the structure of a graded connected cocommutative Hopf algebra by the concatenation product:
$$\mu(y_1...y_n\otimes z_1...z_m)=y_1...y_n\cdot z_1...z_m:=y_1...y_nz_1...z_m$$
and the unshuffling coproduct:
$$\Delta(y_1...y_n):=\sum\limits_{p=0}^n\sum\limits_{I\coprod J=[n]}y_{i_1}...y_{i_p}\otimes y_{j_1}...y_{j_{n-p}},$$
where 
$1\leq i_1<i_2<...<i_p\leq n$, $1\leq j_1<j_2<...<j_{n-p}\leq n$ and $I=\{i_1,...,i_p\}$, $J=\{j_1,...,j_{n-p}\}$.
The antipode is given by: $S(y_1...y_n)=(-1)^ny_n...y_1$, where $y_i\in X$. We refer e.g. to \cite{reutenauer} for further details on the subject and general properties of the tensor algebra viewed as a Hopf algebra and recall simply that the antipode is the convolution inverse of the identity of $T(X)$: $S\ast Id=Id\ast S=\varepsilon$, where $\varepsilon$ is the projection on the scalars $T_0(X)$ orthogonally to the higher degree components $T_n(X),\ n>0$, and the convolution product $f\ast g$ of two linear endomorphisms of $T(X)$ (and, more generally of a Hopf algebra with coproduct $\Delta$ and product $\mu$) is given by: $f\ast g:=\mu\circ (f\otimes g)\circ\Delta$. In particular, $S\ast Id$ is the null map on $T_n(X)$ for $n>0$.

We write $\delta$ an arbitrary derivation of $T(X)$ (in particular $\delta$ acts as the null application on the scalars, $T_0(X)$). The simplest and most common derivations are induced by maps $f$ from $X$ to its linear span: the associated derivation, written $\tilde f$, is then defined by $f(y_1...y_n)=\sum_{i=0}^ny_1...y_{i-1}f(y_i)y_{i+1}...y_n$. 
Since for $a,b\in T(X),\ f([a,b])=[f(a),b]+[a,f(b)]$, where $[a,b]:= ab-ba$, these particular derivations are also derivations of the free Lie algebra over $X$. These are the derivations we will be most interested in in practice (the ones that map the free Lie algebra over $X$ to itself), we call them \it Lie derivations\rm .

In the particular case $f=Id$, we also write $Y$ for $\tilde{Id}$: $Y$ is the graduation operator, $Y(y_1...y_n)=n\cdot y_1...y_n$. When $f=\delta_{x_i}$ ($f(x_i)=x_i, \ f(x_j)=0$ for $j\not= i$), $\tilde f$ counts the multiplicity of the letter $x_i$ in words and is the noncommutative analog of the derivative with respect to $x_i$ of a monomial in the letters in $X$.

\begin{proposition}For arbitrary letters $y_1,...,y_n$ of $X$, we have, for $D_{\delta}:=S\ast \delta$:
 $$D_\delta (y_1...y_n)=[...[[\delta(y_1),y_2],y_3]...,y_n]$$
\end{proposition}

Let us assume, by induction, that the identity holds in degrees $<n$. Then, with the same notation as the ones used to define the coproduct $\Delta$:
$$S\ast \delta (y_1...y_n) = \sum\limits_{p=0}^n(-1)^p y_{i_p}...y_{i_1}\delta(y_{j_1}...y_{j_{n-p}})$$
where $I\coprod J=[n]$. 
Notice that, if $i_p\not=n$, $j_{n-p}=n$ (and conversely).
Therefore, since $\delta$ is a derivation:
$$S\ast \delta (y_1...y_n) = \sum\limits_{p=0}^{n-1}(-1)^p (y_{i_p}...y_{i_1}\delta(y_{j_1}...y_{j_{n-p-1}})y_n$$
$$+ y_{i_p}...y_{i_1}y_{j_1}...y_{j_{n-p-1}}\delta(y_n))$$
$$+\sum\limits_{p=0}^{n-1}(-1)^{p+1}y_ny_{i_p}...y_{i_1}\delta(y_{j_1}...y_{j_{n-p-1}}),$$
where the sums run over the partitions $I\coprod J=[n-1]$.

The first and third term of the summation sum up to $[S\ast \delta (y_1...y_{n-1}),y_n]$ which is, by induction, equal to
$[...[[\delta(y_1),y_2],y_3]...,y_n]$.
The second computes $S\ast Id(y_1...y_{n-1}) \delta(y_n)$, which is equal to $0$ for $n>1$. The Proposition follows.

\begin{corollary}For a Lie derivation $\delta$,
 the map $D_\delta$ maps $T(X)$ to $Lie(X)$, the free Lie algebra over $X$.
 Moreover, we have, for $l\in Lie(X)$, $$D_\delta(l)=\delta(l).$$
\end{corollary}

The first part of the corollary follows from the previous proposition. 
To prove the second part, recall that $l\in Lie(X)$ if and only if $\Delta(l)=l\otimes 1+1\otimes l$. Notice furthermore that, since $D_\delta=S\ast \delta$, we have $\delta=Id\ast D_\delta$. Therefore:
$$\delta(l)=(Id\ast D_\delta)(l)=D_\delta(1)\cdot l+D_\delta(l).$$
The proof follows since $D_\delta(1)=0$.

We recover in particular the theorem of Dynkin \cite{dynkin}, Specht \cite{specht}, Wever \cite{wever} (case $f=Id$) and obtain an extension thereof to the case $f=\delta_{x_i}$. We let the reader derive similar results for other families of Lie derivations.

\begin{corollary}\label{classical}
 We have, for the classical Dynkin operator $D=S\ast Y$ and an arbitrary element $l$ in $T_n(X)\cap Lie(X)$:
 $$D(l)=n\cdot l.$$
 In particular, the operator $\frac D n$ is a projection from $T_n(X)$ onto $T_n(X)\cap Lie(X)$.
\end{corollary}
The definition $D=S*Y$ of the Dynkin operator seems due to von Waldenfels, see \cite{reutenauer}.

Let us write $T_n^{i}(X)$ for the linear span of words over $X$ such that the letter $x_i$ appears exactly $n$ times. The derivation $\tilde \delta_{x_i}$ acts as the multiplication by $n$ on $T_n^{i}(X)$.
\begin{corollary}
 We have for $D_{x_i}:=S\ast \tilde\delta_{x_i}$ and an arbitrary element $l$ in $T_n^i(X)\cap Lie(X)$:
 $$D_{x_i}(l)=n\cdot l.$$
 In particular, the operator $\frac D n$ is a projection from $T_n^i(X)$ onto $T_n^i(X)\cap Lie(X)$.
\end{corollary}

\section{Abstract logarithmic derivatives}

Quite often, the logarithmic derivatives one is interested in arise from dynamical systems and geometry.
We explain briefly why on a fundamental example, the classification of singular vector fields (section~\ref{exemple}).
Although we settle our later computations in the general framework of Lie and enveloping algebras, the reader may keep that motivation in mind. 

The second section (\ref{henv}) shows briefly how to extend the results on generalized Dynkin operators obtained previously in the tensor algebra to the general setting of enveloping algebras. 

We show at last (section~(\ref{idc}) how these results connect to the theory of Rota Baxter algebras, which is known to be the right framework to investigate the formal properties of derivations. Indeed, as we will recall below, Rota-Baxter algebra structures show up naturally when derivations have to be inverted. See also \cite{EGP,EMP,BCEP} for further details on the subject of Rota-Baxter algebras and their applications.

\subsection{An example from the theory of dynamical systems}\label{exemple}

Derivations on graded complete Lie algebras appear naturally in the framework of dynamical systems, especially when dealing with the formal classification (up to formal change of coordinates) of singular vector fields. 

The reader can refer to \cite{Ilya} for an overview and further details on the objects we consider (such as identity-tangent diffeomorphisms or substitution automorphisms) -let us also mention that the reader who is interested only in formal aspects of logarithmic derivatives may skip that section.

A formal singular vector field in dimension $\nu$ is an operator:
\[
X=\sum_{i=1}^{\nu} f_i(x_1,\dots,x_{\nu})\frac{\partial}{\partial x_i}
\]
such that $f_i(0)=0$ for all $i$ (that is $f_i\in \mathbb{C}_{\geq 1}[[x_1,\dots,x_{\nu}]]$). Such operators act on the algebra of formal series in $\nu$ variables. In practice, a vector field  is given by a series of operators such as $x_1^{n_1} \ldots
x_{\nu}^{n_{\nu}} \frac{\partial}{\partial x_i}$ with $n_1 + \ldots + n_{\nu}
\geqslant 1$  that acts on
monomials  $x_1^{m_1} \ldots x_{\nu}^{m_{\nu}}$:
\[ \left( \left. x_1^{n_1} \ldots x_{\nu}^{n_{\nu}} \frac{\partial}{\partial
   x_i} \right) . x_1^{m_1} \ldots x_{\nu}^{m_{\nu}} = m_i\cdot x_1^{m_1 + n_1} \ldots
   x_i^{m_i + n_i - 1} \ldots x_{\nu}^{m_{\nu} + n_{\nu}} \right. \]
so that the total degree goes from $m_1 + \ldots + m_{\nu}$ to $m_1 + \ldots +
m_{\nu} + n_1 + \ldots + n_{\nu} - 1$ and the graduation for such an operator
is then $n_1 + \ldots + n_{\nu} - 1$. 

The $0$ graded component of a vector field $X$ is called the linear part since it can be written $X_0=\sum A_{ij}x_i\frac{\partial}{\partial x_j}$ and a fundamental question in dynamical systems is to decide if $X$ is conjugate, up to a change of coordinates, to its linear part $X_0$.
Notice that:
\begin{itemize}
\item The vector space $L$ of vector fields without linear part (or without component of graduation 0) is a graded complete Lie algebra. %as well as the vector space $L_0$ of linear vector fields.
\item The exponential of a vector field in $L$ gives a one to one correspondence between vector fields and substitution automorphisms on formal power series, that is operators $F$ such that
\[ F(A(x_1,\dots,x_\nu))=A(F(x_1),\dots,F(x_{\nu}))
\]
 where $(F(x_1),\dots,F(x_{\nu}))$ is a formal identity-tangent diffeomorphism.
 \item The previous equation also determines an isomorphism between the Lie group of $L$ and the group of  formal identity-tangent diffeomorphism $G_1$.
 \end{itemize}  
These are essentially the framework (the one of graded complete Lie algebras) and the objects (elements of the corresponding formal Lie groups) that we will consider and investigate in our forthcoming developments.
 
Consider now a vector field $X=X_0 +Y\in L_0\oplus L$ and suppose that it can be linearized by a change of coordinates in $G_1$, or rather by a substitution automorphism $F$ in the Lie group of $L$. It is a matter of fact (see \cite{Ilya}, \cite{SNAG}) to check that the corresponding conjugacy equation reads:
$$
X_0 F=F(X_0 +Y)  \Longleftrightarrow [X_0,F]=FY \Longleftrightarrow ad_{X_0}(F)=FY
$$
This equation, called the homological equation, delivers a derivation $\delta=ad_{X_0}$ on $L$ that is compatible with the graduation. The linearization problem is then obviously related to the inversion of the logarithmic derivation $D_{\delta}(F):=F^{-1}\delta(F)$. 

In the framework of dynamical systems, the forthcoming theorem \ref{invert}  ensures  that if the derivation $ad_{X_0}$ is invertible on $L$, any vector field  $X_0+Y$ can be linearized. This is the kind of problems that can be addressed using the general theory of logarithmic derivatives to be developed in the next sections.

\subsection{Hopf and enveloping algebras}\label{henv}

We use freely in this section the results in \cite{patreu2002} to which we refer for further details and proofs. The purpose of this section is to extend the results in \cite{patreu2002} on the Dynkin operator to more general logarithmic derivatives.

Let $L=\bigoplus\limits_{n\in\NM^\ast}L_n$ be a graded Lie algebra, $\hat L=\prod\limits_{n\in\NM^\ast}L_n$ its completion, $U(L)=\bigoplus\limits_{n\in\NM^\ast}U(L)_n$ the (graded) enveloping algebra of $L$ and $\hat U(L)=\prod\limits_{n\in\NM^\ast}U(L)_n$ the completion of $U(L)$ with respect to the graduation. 

The ground field is chosen to be $\QM$ (but the results in the article would hold for an arbitrary ground field of characteristic zero and, due to the Cartier-Milnor-Moore theorem \cite{mm,patras}, for arbitrary graded connected cocommutative Hopf algebras).

The enveloping algebra $U(L)$ is naturally provided with the structure of a Hopf algebra.
We denote by 
$\e: \QM=U(L)_0\ra U(L)$ the unit of $U(L)$, by $\eta: U(L)\ra \QM$ 
the counit, by $\Delta :U(L)\ra U(L)\ot U(L)$ the coproduct and by 
$\mu :U(L)\ot U(L) \ra U(L)$ the product. 
An element $l$ of $U(L)$ is {\em primitive} if $\Delta (l)=l\ot 1+1\ot l$; the set of
primitive elements identifies canonically with $L$.
Recall that the {\em convolution} product $*$ of linear endomorphisms of $U(L)$ is
defined by $f*g=\mu \circ (f\ot g) \circ \Delta$; $\n:=\e\circ\eta$ is the neutral element of $\ast$. The antipode is written $S$, as usual.

\begin{definition}
An element $f$ of $End(U(L))$ admits $F\in End(U(L))\ot End(U(L))$ as a {\em pseudo-coproduct} if $F \circ \Delta = 
\Delta \circ f$. 
If $f$ admits the pseudo-coproduct $f\ot \nu +\nu\ot f$, we say that $f$ is 
{\em pseudo-primitive}.
\end{definition}

In general, an element of $End(U(L))$ may admit several pseudo-coproducts.
However, this concept is very flexible, as shows the following result \cite[Thm. 2]{patreu2002}.

\begin{proposition}
If $f,g$ admit the pseudo-coproducts $F,G$ and $\a\in \F$, 
then $f+g, \a f, f*g, f\circ g$ admit respectively the pseudo-coproducts
$F+G, \a F, F*G, F\circ G$, where the products $*$ and $\circ$ are naturally extended to $End(U(L))\ot End(U(L))$.

An element $f \in End(U(L))$ takes values in $Prim(U(L))$ if and only if it is pseudo-primitive.
\end{proposition}

Let $\delta$ be an arbitrary derivation of $L$ ($\forall l,l'\in L, \ \delta([l,l'])=[\delta(l),l']+[l,\delta(l')]$). We also write $\delta$ for its unique extension to a derivation of $U(L)$ and write $D_\delta:=S\ast \delta$.
For an element $l\in L$, $\exp(l)$ is group-like ($\Delta(\exp(l))=\exp(l)\otimes \exp(l)$), from which it follows that:
$$D_\delta(\exp(l))=S(\exp(l))\delta(\exp(l))=\exp(-l)\delta(\exp(l)),$$
the (noncommutative) logarithmic derivative of $\exp(l)$ with respect to $\delta$.
We call therefore $D_\delta$ the \it logarithmic derivative \rm of $\delta$.

\begin{proposition}
 The logarithmic derivative $D_\delta$ is a pseudo-primitive: it maps $U(L)$ to $L$.
\end{proposition}

Indeed, $S\otimes S$ is a pseudo-coproduct for $S$ (see \cite{patreu2002}). 
On the other hand $U(L)$ is spanned by products $l_1...l_n$ of elements of $L$. Since $\delta$ is a derivation, we get:
$$\Delta\circ\delta (l_1...l_n)=\Delta(\sum_{i=1}^nl_1...\delta(l_i)...l_n)=(\delta\otimes Id+Id\otimes\delta)\circ\Delta(l_1...l_n),$$
where the last identity follows directly from the fact that the $l_i$ are primitive, which implies that $\Delta(l_1...l_n)$ can be computed by the same formula as the one for the coproduct in the tensor algebra.
In particular, $\delta\otimes Id+Id\otimes\delta$ is a coproduct for $\delta$.
We get finally:
$$\Delta\circ D_\delta=\Delta\circ(S\ast\delta)=(S\otimes S)\ast(\delta\otimes Id+Id\otimes\delta)\circ\Delta$$
$$=(D_\delta\otimes \n+\n\otimes D_\delta)\circ\Delta ,$$
from which the proof follows.

\begin{proposition}\label{primit}
For $l\in L$, we have $\delta(l)=D_\delta(l)$. In particular, when $\delta$ is invertible on $L$, $D_\delta(l)$ is a projection from $U(L)$ onto $L$.
\end{proposition}

The proof is similar to the one in the free Lie algebra.
We have $D_\delta(l)=(S*\delta)(l)=\pi\circ (S\ot \delta)\circ \Delta (l)=\pi\circ (S\ot \delta)(l\ot 1+1\ot l)$
$=\pi (S(l)\ot \delta(1)+S(1)\ot \delta(l))=\delta(l)$, since $\delta(1)=0$ and $S(1)=1$.

\subsection{Integro-differential calculus}\label{idc}

The notation are the same as in the previous section, but we assume now that the derivation $\delta$ is invertible on $L$ and extends to an invertible derivation on $U(L)^+:=\bigoplus\limits_{n\geq 1}U(L)_n$. The simplest example is provided by the graduation operator $Y(l)=n\cdot l$ for $l\in L_n$ (resp. $Y(x)=n\cdot x$ for $x\in U(L)_n$). This includes the particular case, generic for various applications to the study of dynamical systems, where $L$ is the graded Lie algebra of polynomial vector fields spanned linearly by the $x^n\partial_x$ and $\delta :=x\partial_x$ acting on $P(x)\partial_x$ as $\delta(P(x)\partial_x):=xP'(x)\partial_x$.

We are interested in the linear differential equation
\begin{equation}\label{eqdif}
 \delta \phi=\phi\cdot x,\ x\in L,\ \phi\in 1+U(L)^+.
\end{equation}
The inverse of $\delta$ is written $R$ and satisfies, on $U(L)^+$, the identity:
$$R(x)R(y)=R(R(x)y)+R(xR(y)),$$
that follows immediately from the fact that $\delta$ is an invertible derivation. In other terms, $U(L)^+$ is provided by $R$ with the structure of a weight $0$ Rota-Baxter algebra and solving (\ref{eqdif}) amounts to solve the so-called Atkinson recursion:
$$\phi =1+R(\phi\cdot x).$$
We refer to \cite{EGP,EMP} for a detailed study of the solutions to the Atkinson recursion and further references on the subject.
Perturbatively, the solution is given by the Picard (or Chen, or Dyson... the name given to the expansion depending on the application field) series:
$$\phi = 1+\sum\limits_{n\geq 1}R^{[n]}(x),$$
where $R^{[1]}(x)=R(x)$ and $R^{[n]}(x):=R(R^{[n-1]}(x)x)$.

Since the restriction to the weight $0$ is not necessary for our forthcoming computations, we restate the problem in a more general setting and assume from now on that $R$ is a weight $\theta$ Rota-Baxter (RB) map on $U(L)^+$, the enveloping algebra of a graded Lie algebra. That is, we assume that:
$$R(x)R(y)=R(R(x)y)+R(xR(y))-\theta R(xy).$$
This assumption allows to extend vastly the scope of our forthcoming results since the setting of Rota-Baxter algebras of arbitrary weight includes, among others, renormalization in perturbative quantum field theory and difference calculus, the later setting being relevant to the study of diffeomorphisms in the field of dynamical systems. We refer in particular to the various works of K. Ebrahimi-Fard on the subject (see e.g. \cite{EMP,BCEP} for various examples of RB structures and further references).
We assume furthermore that $R$ respects the graduation and restricts to a linear endomorphism of $L$.

\begin{lemma}
 The solution to the Atkinson recursion is a group-like element in $1+U(L)^+$. In particular, $S(\phi)=\phi^{-1}$.
\end{lemma}

Recall from \cite{patreu2002} and \cite{EFGP} that the generalized Dynkin operator $D:=S\ast Y$ (the convolution of the antipode with the graduation map in $U(L)$) maps $U(L)$ to $L$ and, more specifically, defines a bijection between the set of group-like elements in $\hat U(L)$ and $\hat L$.
The inverse is given explicitly by (\cite[Thm. 4.1]{EFGP}):
$$D^{-1}(l)=1+\sum_{n\in\NM^\ast}\sum_{k_1+...+k_l=n}\frac{l_{k_1}\cdot ...\cdot l_{k_l}}{k_1(k_1+k_2)...(k_1+...+k_l)},$$
where $l_n$ is the component of $l\in L$ in $L_n$.
According to \cite[Thm. 4.3]{EGP}, when $l=D(1+\sum\limits_{n\geq 1}R^{[n]}(x))$ we also have:
$$1+\sum\limits_{n\geq 1}R^{[n]}(x)=1+\sum_{n\in\NM^\ast}\sum_{k_1+...+k_l=n}\frac{l_{k_1}\cdot ...\cdot l_{k_l}}{k_1(k_1+k_2)...(k_1+...+k_l)},$$
that is, since $l\in L$ by Prop.~\ref{primit}, $1+\sum\limits_{n\geq 1}R^{[n]}(x)$ is a group-like element in $U(L)$. The last part of the Lemma is a general property of the antipode acting on a group-like element; the Lemma follows.

We are interested now in the situation where another Lie derivation $d$ acts on $U(L)$ and commutes with $R$ (or equivalently with $\delta$ when $R$ is the inverse of a derivation). A typical situation is given by Schwarz commuting rules between two different differential operators associated to two independent variables.

\begin{theorem}
 Let $d$ be a graded derivation on $U(L)$ commuting with the weight $\theta$ Rota-Baxter operator $R$. Then, for $\phi$ a solution of the Atkinson recursion as above, we have:
 $$D_d(\phi)=\phi^{-1}\cdot d(\phi)=\sum\limits_{n\geq 1}R_d^{[n]}(x)$$
 with $I_d^{[1]}=d(x)$, $R_d^{[1]}(x)=R(d(x))$, $I_d^{[n+1]}(x)=[R_d^{[n]}(x),x]+\theta x\cdot I_d^{[n]}(x)$ and $R_d^{[n+1]}(x)=R(I_d^{[n+1]}(x))$.
\end{theorem}

Notice, although we won't make use of this property, that the operation $x\circ y:=[R(x),y]+\theta y\cdot x$ showing up implicitly in this recursion is a preLie product, see e.g. \cite{EGP}. In particular, for  $y$ the solution to the preLie recursion:
$$y=d(x)+y\circ x,$$
we have: $D_d(\phi)=R(y)$.

The first identity $D_d(\phi)=\phi^{-1}\cdot d(\phi)$ follows from the definition of the logarithmic derivative $D_d:=S\ast d$ and from the previous Lemma.

The second is equivalent to
$d(\phi)=\phi\cdot \sum\limits_{n\geq 1}R_d^{[n]}(x)$, that is:
$$d(R^{[n]}(x))=R_d^{[n]}(x)+\sum_{i=1}^{n-1}R^{[i]}(x)R_d^{[n-i]}(x).$$
For $n=1$, the equation reads $d(R(x))=R(d(x))$ and expresses the commutation of $d$ and $R$. 

The general case follows by induction.
Let us assume that the identity holds for the components in degree $n<p$ of $D_d(\phi)$. We summarize in a technical Lemma the main ingredient of the proof. Notice that the Lemma follows directly from the Rota-Baxter relation and the definition of $R_d^{[n+1]}(x)$.

\begin{lemma}
 We have, for $n,m\geq 1$:
 $$R(R^{[m]}(x)\cdot ([R_d^{[n]}(x),x]+\theta x\cdot I_d^{[n]}(x)))=R^{[m]}(x)R_d^{[n+1]}(x)$$
 $$-R(R^{[m-1]}(x)\cdot x\cdot R_d^{[n+1]}(x))+\theta R(R^{[m-1]}(x)\cdot x\cdot ([R_d^{[n]}(x),x]+\theta x\cdot I_d^{[n]}(x)).$$
\end{lemma}

We have, for the degree $p$ component of $D_d(\phi)$, using the Lemma to rewrite $R(R^{[m]}(x)\cdot R_d^{[n]}(x)\cdot x)$:
$$d(R^{[p]}(x))=d(R(R^{[p-1]}(x)\cdot x))=R((dR^{[p-1]}(x))\cdot x+R^{[p-1]}(x)\cdot dx)$$
$$=R(R^{[p-1]}(x)\cdot dx)+R((R_d^{[p-1]}(x)+\sum_{k=1}^{p-2}R^{[p-1-k]}(x)\cdot R_d^{[k]}(x)   )\cdot x)$$
$$=R(R^{[p-1]}(x)\cdot dx)+R(R_d^{[p-1]}(x)\cdot x)+\sum\limits_{k=1}^{p-2}[ R^{[p-1-k]}(x)\cdot R_d^{[k+1]}(x)$$
$$-R(R^{[p-2-k]}(x)\cdot x\cdot R_d^{[k+1]}(x))+\theta R(R^{[p-2-k]}(x)\cdot x\cdot ([R_d^{[k]}(x),x]+\theta x\cdot I_d^{[k]}(x))$$
$$+R(R^{[p-1-k]}(x)\cdot x\cdot R_d^{[k]}(x))-\theta R(R^{[p-1-k]}(x)\cdot x\cdot I_d^{[k]}(x))].$$
The fourth and sixth terms cancel partially and add up to $R(R^{[p-2]}(x)\cdot x\cdot R_d^{[1]}(x))-R(x\cdot R_d^{[p-1]}(x))$. The fifth and last terms cancel partially and add up to $\theta R(x\cdot I_d^{[p-1]}(x))-\theta R(R^{[p-2]}(x)\cdot x\cdot I_d^{[1]}(x))$.
In the end, we get:
$$d(R^{[p]}(x))=\sum\limits_{k=1}^{p-2}R^{[p-1-k]}(x)\cdot R_d^{[k+1]}(x)+[R(R^{[p-2]}(x)\cdot x\cdot R_d^{[1]}(x))+$$
$$R(R^{[p-1]}(x)\cdot dx)-\theta R(R^{[p-2]}(x)\cdot x\cdot I_d^{[1]}(x)]+[R(R_d^{[p-1]}(x)\cdot x)-R(x\cdot R_d^{[p-1]}(x))$$
$$+\theta R(x\cdot I_d^{[p-1]}(x))]
$$
Using the RB identity for the expressions inside brackets, we get finally:
$$d(R^{[p]}(x))=\sum\limits_{k=1}^{p-2}R^{[p-1-k]}(x)\cdot R_d^{[k+1]}(x)+R^{[p-1]}(x)R_d^{[1]}(x)+R_d^{[p]}(x),$$
from which the Theorem follows.

\section{Magnus-type formulas}

The classical Magnus formula relates ``logarithms and logarithmic derivatives'' in the framework of linear differential equations. That is, it relates explicitly the logarithm $\log(X(t))=:\Omega(t)$ of the
solution to an arbitrary matrix (or operator) differential equation 
$X'(t)=A(t)X(t),\ X(0)=1$
to the infinitesimal generator $A(t)$ of the differential equation:
$$\Omega'(t)=\frac{ad_{\Omega(t)}}{\exp^{ad_{\Omega(t)}}-1}A(t)=A(t)+\sum\limits_{n>0}\frac{B_n}{n!}ad_{\Omega(t)}^n(A(t)),$$
where $ad$ stands for the adjoint representation and the $B_n$ for the Bernoulli numbers.

The Magnus formula is a useful tool for numerical applications (computing the logarithm of the solution improves the convergence at a given order of approximation). It has recently been investigated and generalized from various points of view, see e.g. \cite{EM2}, where the formula has been extended to general dendriform algebras (i.e. noncommutative shuffle algebras such as the algebras of iterated integrals of operators), \cite{CP} where the algebraic structure of the equation was investigated from a purely preLie algebras point of view, or \cite{BCEP} where a generalization of the formula has been introduced to model the commutation of time-ordered products with time-derivations.

The link with preLie algebras follows from the observation that (under the hypothesis that the integrals and derivatives are well-defined), for arbitrary time-dependent operators, the preLie product
$$M(t)\curvearrowleft N(t):=\int_0^t [N(u),M'(u)]du$$
satisfies $(M(t)\curvearrowleft N(t))'=ad_{N(t)}M'(t)$. The Magnus formula rewrites therefore:
$$\Omega'(t)=\left(\int\limits_0^tA(x)dx\curvearrowleft\left( \frac{\Omega}{\exp(\Omega)-1} \right)\right)'$$
where $\frac{\Omega}{\exp(\Omega)-1}$ is computed in the enveloping algebra of the preLie algebra of time-dependent operators. 

Here, we would like to go one step further and extend the formula to general logarithmic derivatives in the suitable framework in view of applications to dynamical systems and geometry, that is, the framework of enveloping algebras of graded Lie algebras and Lie derivations actions.
Notations are as before, that is $L$ is a graded Lie algebra and $\delta$ a graded Lie derivation (notice that we do not assume its invertibility on $L$ or $U(L)$).

\begin{lemma}
 For $l\in L$, $k\geq 1$ and $x\in U(L)$, we have:
 $$x\cdot l^k=\sum\limits_{i=0}^k{{k}\choose{ i}}l^{k-i}\cdot (-ad_{l})^i(x). $$
\end{lemma}

The proof is by induction on $n$. The formula holds for $n=1$: $xl=-[l,x]+lx$. Let us assume that it holds for an arbitrary $n<p$. Then we have:
$$x\cdot l^p=  (x\cdot l^{p-1})\cdot l= (\sum\limits_{i=0}^{p-1}{{p-1}\choose{ i}}l^{p-1-i}\cdot (-ad_{l})^i(x))\cdot l$$
$$=\sum\limits_{i=0}^{p-1}{{p-1}\choose{ i}}l^{p-1-i}\cdot (-ad_{l})^{i+1}(x) +\sum\limits_{i=0}^{p-1}{{p-1}\choose{ i}} l^{p-i}\cdot(-ad_{l})^{i}(x).$$
The identity follows then from Pascal's triangular computation of the binomial coefficients.

\begin{theorem}\label{magnus}
 For $l\in L$, we have:
 $$D_\delta (\exp(l))=\frac{\exp(-ad_l)-1}{-ad_l}\delta(l).$$
\end{theorem}

Indeed, from the previous formula we get, 
$$d(\exp(l))=\sum\limits_{n\geq 1}\frac{1}{n!}d(l^n)= \sum\limits_{n\geq 1}\frac{1}{n!}\sum_{k=0}^{n-1}l^{n-1-k}d(l)l^k$$
$$=\sum\limits_{n\geq 1}\sum\limits_{i=0}^{n-1}\frac{1}{n!}(\sum\limits_{k=i}^{n-1}{k\choose i})l^{n-1-i}(-ad_l)^i(d(l))$$
$$=\sum\limits_{n\geq 1}\sum\limits_{i=0}^{n-1}\frac{1}{n!}{n\choose i+1}l^{n-1-i}(-ad_l)^i(d(l))$$
$$=\sum\limits_{i\geq 0}\sum\limits_{n\geq i+1}\frac{1}{(i+1)!(n-1-i)!}l^{n-1-i}(-ad_l)^i(d(l))$$
$$=\exp(l)(\sum\limits_{i\geq 0}\frac{1}{(i+1)!}(-ad_l)^i)d(v).$$
Since $\exp(l)$, being the exponential of a Lie element is group-like, $D_\delta(\exp(l))= \exp({-l})d(\exp(l))$ and
the theorem follows.

Let us show as a direct application of the Thm (\ref{magnus}) how to recover the classical Magnus theorem (other applications to mould calculus and to the formal and analytic classification of vector fields are postponed to later works).

\begin{example}
 Let us consider once again an operator-valued linear differential equation
 $X'(t)=X(t)\lambda A(t),\ X(0)=1$. Notice that the generator $A(t)$ is written to the right, for consistency with our definition of logarithmic derivatives $D_\delta =S\ast d$. All our results can of course be easily adapted to the case $d\ast S$ (in the case of linear differential equations this amounts to consider instead $X'(t)=A(t)X(t)$), this easy task is left to the reader. Notice also that we introduce an extra parameter $\lambda$, so that the perturbative expansion of $X(t)$ is a formal power series in $\lambda$.
 
Consider then the Lie algebra $O$ of operators $M(t)$ (equipped with the bracket of operators) and the graded Lie algebra $L=\bigoplus_{n\in \NM^\ast}\lambda^n O=O\otimes \lambda \CM[\lambda]$. Applying the Thm (\ref{magnus}), we recover the classical Magnus formula.
\end{example}

Although a direct consequence of the previous theorem (recall that any group-like element in the enveloping algebra $U(L)$ can be written as an exponential), the following proposition is important and we state it also as a theorem.

\begin{theorem} \label{invert}
 When $\delta$ is invertible on $L$, the logarithmic derivative $D_\delta$ is a bijection between the set of group-like elements in $U(L)$ and $L$.
\end{theorem}

Indeed, for $h\in L$ the Magnus-type equation in $L$
$$l=\delta^{-1}(\frac{-ad_l}{\exp (-ad_l)-1}(h))$$
has a unique recursive solution 
$ l\in L$ such that  $exp(l)=D_\delta^{-1}(h)$.

% 
% dans une algèbre de Hopf cocommutative connexe (moules, etc), on aurait pour x primitif, D dérivation et \phi un automorphisme d'algèbre
% 
% (D\ast (\phi\circ S))(exp(x))=\frac{exp(ad_{\phi(x)}-1)}{ad_{\phi(x)}} Dx
% 

%%%%%%%%%%%%%%%%%%%%%%%%%%%%%%%%%%%%%%%%%%%%%%%%%%%%%%%%%%%%%%%%%%%%%%%%%%%%%%%%%%%%%%%%%%%%%%%%%%%%%%%%

\end{document}